\documentclass[11pt]{article}
\usepackage{amsmath,amssymb}

\newcommand{\C}{\mathbb{C}}

\newcommand{\R}{\mathbb{R}}

\newcommand{\bfg}{\mathbf{g}}
\newcommand{\rr}{\mathbf{r}}

\newcommand{\nn}{\mathbf{n}}

\newcommand{\ep}{\varepsilon}

\newcommand{\ph}{\varphi}
\newcommand{\rd}{\partial}
\newcommand{\rdo}{Rad\'{o} }
\newcommand{\rdos}{Rad\'{o}'s }

\numberwithin{equation}{section}

\begin{document}

\title{  The work of Jesse Douglas on Minimal Surface}

\author{Jeremy Gray, \\
Centre for the History of the Mathematical Sciences,\\
Open University, Milton Keynes, MK7 6AA, U.K., \and
Mario Micallef, \\
Mathematics Institute, University of Warwick, Coventry, CV4 7AL,
U.K.}

\maketitle This paper is dedicated to the memory of our student,
Adam Merton (1920-1999) without whose interest in minimal
surfaces and the work of Jesse Douglas we would never have
embarked on this study.

\section{Introduction }
The demonstration of the existence of a least-area surface
spanning a given contour has a long history. Through the
soap-film experiments of Plateau, this existence problem
became known as the Plateau Problem.
By the end of the 19th century, the list of contours for
which the Plateau problem could be solved was intriguingly close
to that considered by Plateau, namely polygonal contours, lines
or rays  extending to infinity and circles. Free boundaries on
planes were also considered. But even these contours were
required to have some special symmetry, for example a pair of
circles had to be in parallel planes. These solutions were
obtained by finding the pair of holomorphic functions (which, for
the experts, we shall refer to as the conformal factor $f$ and
the Gauss map $g$) required in the so-called Weierstrass-Enneper
representation of a minimal surface.  In 1928, Garnier completed
a programme started independently by Riemann and Weierstrass in
the 1860s  and then continued by Darboux towards the end of the
19th Century for finding $f$ and $g$ for an arbitrary polygonal
contour. In his 92-page paper~\cite{garn} (which is very
difficult to read and which, apparently has never been fully
checked) Garnier also employed a limiting process to extend his
solution of the Plateau problem to contours which consisted of a
finite number of unknotted arcs with bounded curvature. Garnier's
work was soon eclipsed by the works of \rdo and Jesse Douglas
about which there is considerable amount of inaccurate
information in the literature and on which we hope to shed some
light in this article. In particular, we challenge the popular
belief that Douglas arrived at his functional for solving the
Plateau Problem by direct consideration of Dirichlet's integral
and its relation to the area functional. Douglas was awarded one of
the first Fields Medals for his work on the Plateau problem.
There are many amusing aspects of the Fields Medal ceremony
at which Douglas was awarded his prize. We simply mention that
Wiener collected the medal on behalf of Douglas,
even though Douglas did attend the International Congress.
And in the address, Carath\'{e}odory described a method
for finding a minimal surface that is due to \rdo
(different from the one sketched below) but gave the impression that
it was due to Douglas! Full historical and
mathematical details will appear in \cite{jm}.

\section{Minimal Surfaces as Conformal Harmonic Maps}

Recall that $S \subset \R^n$ is called a regular (immersed)
surface if there exists a surface (a manifold of real dimension
2) $\Sigma$, possibly with boundary, and a map $\rr \colon \Sigma
\to \R^n$ of class $C^1$ up to the boundary such that $\rr(\Sigma)
= S$ and the differential of $\rr$ has maximal rank (and so, the
induced map on tangent spaces is an injection).\footnote{Not much
is lost for the present purposes if one restricts attention to
the simple case when $\Sigma$ is a domain in $\R^2$.} The map
$\rr$ is called a parameterisation of $S$. The first fundamental
form, $I$, of $\rr$ with respect to local coordinates $(u,v)$ on
$\Sigma$ is the symmetric matrix defined by:
\begin{equation}\label{efg}
I := \begin{pmatrix}
E & F \\
F & G
\end{pmatrix}
\quad\text{where}\quad E := \frac{\rd \rr}{\rd u} \cdot \frac{\rd
\rr}{\rd u}, \ F := \frac{\rd \rr}{\rd u} \cdot \frac{\rd
\rr}{\rd v}, \ G := \frac{\rd \rr}{\rd v} \cdot \frac{\rd
\rr}{\rd v}.
\end{equation}
The condition that $\rr$ is of maximal rank is equivalent to the
positive definiteness of $I$. The area of $S$ is defined by the
integral
\begin{equation}\label{area}
\text{Area}(S) := \iint_{\Sigma} \sqrt{g}\,dudv \quad\text{where}
\quad g := EG - F^2 = \text{det}\,I.
\end{equation}
A \emph{critical point} of this area functional is called a
\emph{minimal} surface. It has zero mean curvature,
i.~e. the shape operator $B$ (presently defined for
a surface $S$ in $\R^3$) has zero trace.
The Gauss map $G \colon S \to S^2$, where $S^2$ is
the round unit sphere in $\R^3$, is simply defined by
$G(p) := \nn(p)$, where $\nn$ is the positively oriented unit normal
of $S$. The tangent space of $S^2$ at $G(p)$ and
the tangent space $T_pS$ of $S$ at $p$ may be identified because
both are orthogonal to $\nn(p)$. Therefore, the differential of
$G$ at $p$ may be viewed as an endomorphism $B(p)$ of $T_p(S)$.
The second fundamental form $II$ and the shape operator $B$
are related via the first fundamental form $I$ by
$II(X,Y) = -I(BX,Y) \ \forall \, X,Y \in T_pS$.

It is important (for the purposes of \S \ref{rmt_lap}) to
realise that a minimal surface need not actually minimise area in
the same way that a geodesic need not minimise length.
For instance, consider a catenoid and a cylinder with a common axis
of rotation and which intesect in a pair of circles. If the radius
of the cylinder is not too large, the portion of the catenoid
within the cylinder will have less area than
the portion of the cylinder bounded by the circles.
But if the radius of the cylinder is sufficiently large,
the opposite will be true, and the area minimiser is then
another catenoid which is `fatter' than the one initially considered.

We shall now assume $S$ to be the image of an immersion $\rr
\colon D \to \R^n$ of the open unit disc $D \subset \R^2 = \C$.
It is a fundamental theorem in surface theory (due to Gauss in
the analytic case) that $S$ can be realised as the image of a
conformal map of $D$ --- so we may take the immersion $\rr$ to be
conformal in $u$ and $v$. This means that $E = G$ and $F = 0$.
The coordinates $(u,v)$ are then called isothermal parameters.
The condition that $S$ be minimal, i.~e. that its mean curvature
vanishes, is equivalent to requiring each component $x_i,\ 1
\leqslant i \leqslant n$, of $\rr$ be a harmonic function of the
isothermal coordinates; the map $\rr$ is then called a harmonic
map (called a harmonic vector in the days of Douglas). If the
condition of regularity of $\rr$ is dropped, i.~e., $\rr \colon D
\to \R^n$ is harmonic and satisfies $E=G$ and $F=0$ but $E$ (and
therefore $G$) is allowed to vanish, then $\rr$ parameterises, in
modern terminology, a \emph{generalised, or branched, minimal
disc}. The (isolated) points at which $E$ vanishes are called
branch points.

After the work of Riemann and Weierstrass, the
Plateau problem acquired the following formulation (with $n=3$):

\medskip
\noindent\textbf{Plateau Problem} \ \textit{Given a contour
$\Gamma \subset \R^n$, span $\Gamma$ by a (possibly) branched
minimal disc, i.~e., find $\rr \colon \overline{D} \to \R^n$
which is continuous on the closed unit disc $\overline{D}$,
harmonic and conformal (away from branch points) on the interior
$D$ and whose restriction to $\rd D$ parameterises $\Gamma$.}

As noted above, only modest progress was made on this problem in
the several decades from the 1860s until the end of the 1920s.
Then two mathematicians working independently and with quite
different methods were able to make decisive progress and produce
general methods for solving it. They were Tibor \rdo in Europe
and Jesse Douglas in the United States.

\section{\rdos solution of the Plateau problem} \label{rplat}
\rdo based his solution of the Plateau problem in \cite{rad30b} on
the uniformisation theorem of Koebe.

Given a polygonal contour $\Gamma \subset \R^3$, define
\[
\lambda := \inf\{\text{Area}(\Pi) \mid \Pi \text{ is a polyhedral
surface spanning $\Gamma$}\}.
\]
Then, for each $\sigma > 0$, there exists a polyhedron $\Pi_{\sigma}$
spanning $\Gamma$ whose area is less than $\lambda + \sigma$.
By the uniformisation theorem of Koebe (see Koebe's~\cite{koe17}
and \cite{koe27}), $\Pi_{\sigma}$ admits an isothermic
parameterisation $\bar{\rr}_{\sigma} \colon \overline{D} \to
\R^3$. Let $\rr_{\sigma}$ be the harmonic extension of
$\bar{\rr}_{\sigma}$ restricted to the unit circle $C$, the
boundary of $\overline{D}$. By a lemma on harmonic surfaces, \rdo
asserted the existence of a polyhedron $\Pi_{\sigma}^*$ whose
area differs from that of $\rr_{\sigma}(\overline{D})$ by no more
than $\sigma$. Let $\sqrt{EG-F^2}$ and
$\sqrt{\bar{E}\bar{G}-\bar{F}^2}$ denote the area elements of
$\rr_{\sigma}$ and $\bar{\rr}_{\sigma}$ respectively. Since the
Dirichlet energy $\frac12 \iint_{\overline{D}}\bar{E} + \bar{G}$
of $\bar{\rr}_{\sigma}$ is greater than or equal to the Dirichlet
energy $\frac12 \iint_{\overline{D}}E + G$ of $\rr_{\sigma}$ we
have the following chain of inequalities:
\begin{multline*}
\lambda + \sigma >
\iint_{\overline{D}}\sqrt{\bar{E}\bar{G}-\bar{F}^2}
= \frac12 \iint_{\overline{D}}\bar{E} + \bar{G} \\
\geqslant \frac12\iint_{\overline{D}}E + G \geqslant
\iint_{\overline{D}}\sqrt{EG-F^2} > \text{Area}(\Pi_{\sigma}^*) -
\sigma \geqslant \lambda - \sigma.
\end{multline*}
Therefore, by taking $\sigma$ sufficiently small, one can find,
for each $\ep >0$, a harmonic vector $\rr_{\ep} \colon D \to \R^3$ which
\begin{enumerate}
\item[(i)]extends continuously to the closed unit
disc $\overline{D}$ so that its restriction to $\rd D$
parameterises $\Gamma$ and
\item[(ii)]is approximately conformal in the sense that
\begin{equation}\label{apsol}
\iint_{\overline{D}} |F| < \ep \qquad\text{and}
\qquad \iint_{\overline{D}} (E^{1/2} - G^{1/2})^2 < \ep.
\end{equation}
\end{enumerate}
Let $(\ep_n)$ be a sequence of positive numbers decreasing to 0 and
denote $\rr_{\ep_n}$ more simply by $\rr_n$. \rdo showed that
a subsequence of $(\rr_n)$ converges, up to
reparameterisation by M\"{o}bius transformations of $D$, to a
generalised minimal surface spanning $\Gamma$. This limiting
argument is one of \rdos major achievements. He had previously
established the following variant of it in \cite{rad30a}.

\medskip
\noindent\textbf{Approximation Theorem} \textit{Let $\Gamma_n$ be
a sequence of simple closed curves of uniformly bounded length,
for each of which the Plateau Problem is solvable. If $\Gamma_n$
converges (in the sense of Fr\'{e}chet) to a simple closed curve
$\Gamma$ then the Plateau problem for $\Gamma$ is solvable.}

\medskip
The proof of this Approximation Theorem was not only
\emph{much} simpler than Garnier's limiting argument but,
together with \cite{rad30b}, it also provided a solution to
Plateau's problem for any rectifiable contour.

\rdos use of polyhedral surfaces is highly reminiscent of
Lebesgue's definition of area of a surface as the infimum,
over all sequences, of $\liminf$ of the areas of a sequence of
polyhedra tending to the surface. (The analogous definition for
curves coincides with the usual limsup definition. If the
surface admits a Lipschitz parameterisaton, or even a
parameterisation which has square integrable first derivatives,
then Lebesgue's definition agrees with the one given by the usual
double integral.) It is surprising that \rdo did not connect
his method of proof with Lebesgue's definition of area
until he gave a colloquium on his result at Harvard.
We refer the reader to \cite{jm} for details of this story.
We simply note here that, once \rdo made this connection,
he wrote \cite{rad30c} in which he showed that
his solution of the Plateau problem has least area
among discs spanning $\Gamma$.

\rdos strategy marked a total departure from
the Riemann-Weierstrass-Darboux programme and
it was strikingly original at the time.
We now turn to the solution of Jesse Douglas,
who also abandoned the Riemann-Weierstrass-Darboux programme
but in a direction different from that of Rad\'{o}.

\section{Douglas's solution of the Plateau problem} \label{Dsol}
Douglas's first full account of his solution of the Plateau
problem appears in his 59-page paper `Solution of the Problem of
Plateau', published in the \textit{Transactions of the American
Mathematical Society} in January 1931. However,
he had made several announcements previously at meetings of
the American Mathematical Society, starting in late 1926.
In these announcements, which have essentially been forgotten,
Douglas made clear that his plan was to find
a parameterisation $\bfg^*$ of the given contour $\Gamma$ so that
the harmonic extension $\rr^*$ of $\bfg^*$ is conformal.
In the famous 1931 paper, Douglas achieved this by
minimising his so-called $A$-functional:
\begin{equation}\label{Adef}
A(\bfg) := \frac{1}{4 \pi}\int_0^{2 \pi} \int_0^{2 \pi}
\frac{\sum_{i=1}^n[g_i(\theta) - g_i(\ph)]^2} {4 \sin^2
\frac{\theta - \ph}{2}} \, d \theta d \ph .
\end{equation}
He gives no clue as to how he thought of this functional but, in
Part III of his paper, he does show that $A(\bfg)$ is equal to
the Dirichlet energy $\frac12\iint_{\overline{D}} (E + G)$ of the
harmonic extension $\rr_{\bfg}$ of $\bfg$. It has always been
assumed that this was Douglas's starting point but his
announcements in the \textit{Bulletin} indicate otherwise.
Moreover, Douglas repeatedly refuted claims made by \rdo and
Courant that Douglas's method essentially amounted to
implementing Dirichlet's principle. It must be remembered that,
in 1931, Dirichlet's principle was not yet firmly established and
Douglas was keen, indeed all too keen, to emphasize that his
$A$-functional, being a 1-dimensional integral which did not
involve any derivatives, did not suffer from all the difficulties
that plagued Dirichlet's integral at the time.
Specifically, as Douglas remarked, the Dirichlet integral could
not be shown to  attain its lower bound, whereas his
$A$-functional, being lower semi-continuous on a sequentially
compact space, necessarily attained its minimum. Details of the
exchanges between \rdo and Douglas and Courant and Douglas can be
found in \cite{jm}.

\subsection{Conformality of a harmonic surface via an integral equation
for the parameterisation of the boundary} \label{conf_inteq} In
the Abstract (1927, 32, 143-4)\footnote{Abstracts will be
referred to by the year of publication and number in the Bulletin
of the American Mathematical Society, and listed separately at
the end of this paper.} Douglas claimed that, if $t \mapsto \bfg(t)
\colon \R \cup \{\infty\} \to \Gamma \subset \R^3$ is a
parameterisation of $\Gamma$ and if $\ph \colon \R \cup
\{\infty\} \to \R \cup \{\infty\}$ is a homeomorphism which
solves the following integral equation,
\begin{equation} \label{hiq}
\int_\Gamma \frac{K(t,\tau)}{\ph(t) - \ph(\tau)} d \tau = 0,
\qquad \text{where $K(t,\tau)=\bfg'(t) \cdot \bfg'(\tau)$}
\end{equation}
then the harmonic extension of $\bfg \circ \ph^{-1}$ to the
upper half plane defined by means of Poisson's integral
is the required conformal harmonic representation of
a minimal surface spanning $\Gamma$. (Douglas was actually
less precise than this in the abstract but it is clear that
this is what he meant.) If $\Gamma$ is parameterised by
a map from the unit circle $C$ (and a surface $S$ spanning $\Gamma$
is parameterised by a map from the disc) then
the integral equation \eqref{hiq} becomes
\begin{equation} \label{kinteq}
\int_0^{2 \pi} K(t, \tau) \cot(\tfrac12\ph(t) - \tfrac12\ph(\tau))
\,d \tau = 0 .
\end{equation}
Douglas never published a proof of his claim.
A derivation of \eqref{kinteq} can be found in \cite{jm}.

Integral equations were a highly active topic of research at the
time. Therefore, the formulation of the conformality requirement
on the parameterisation of the minimal surface in terms of an
integral equation would have been perceived as valuable, even
though the integral equation itself is intractable except for
special contours. More importantly, it was a \emph{new,
non-algebraic} formulation of the conformality condition.
As already stated, seeking the necessary special parameterisation of
the contour as the solution of an integral equation was
a significant departure from the then current strategy,
adopted by Garnier, of seeking Weierstrass-Enneper data from
the contour.

\subsection{The integral equation as
the Euler-Lagrange equation of a functional} \label{inteq_el}
Douglas was stuck for a while on how to solve \eqref{kinteq}
for a general contour. His important breakthrough,
which he announced in his~\cite{dg28a} published in
the July-August 1928 issue of the \emph{Bulletin of
the American Mathematical Society}, came when he realised that
\eqref{kinteq} is the Euler-Lagrange equation of
the first version of his later-to-be-famous $A$-functional:
\begin{equation} \label{aph}
A(\ph) := -\int_0^{2\pi}\int_0^{2\pi} K(t,\tau) \log \, \sin
\tfrac12|\ph(t) - \ph(\tau)| \, dt \, d \tau \,.
\end{equation}
Furthermore, he stated that he could use Fr\'{e}chet's
compactness theory of curves to assert the existence of
a minimizer $\ph^*$ of the $A$-functional, at least in the case
that $K(t,\tau)$ is positive for all values of $t$ and
$\tau$. This positivity requirement on $K$ rendered
this variational problem inapplicable to the Plateau problem but
he overcame this problem when he discovered that he could employ
the functional \eqref{Adef} instead of \eqref{aph}.

Douglas gave many talks on his solution of the Plateau problem,
especially in 1929 during his European tour.
A 4-page document written in April 1929 seems to be
his notes for one of these talks. It is contained in a black notebook,
now in City College New York, headed Lebesgue 4i\`{e}me conference,
13 Dec 1928. The relevant pages are between material dated
11 April 1929 and 29 April 1929 and the first page is reproduced
on the front cover. The significance of seeking the minimal surface
in any number of dimensions and of enriching the problem of Plateau
by requiring that the minimal surface be a conformal image of a domain
is explained in the next subsection. On the second page, he wrote:
``This additional requirement is a \underline{natural} complement of
the usual enunciation which demands only $M$, the minimal surface,
in the  sense that the problem is really simpler with it,
than without  it" [emphasis in original]. The document ends with
the integral equation \eqref{kinteq}. Presumably,
he would have gone on to state that \eqref{kinteq} is
the Euler-Lagrange equation of \eqref{aph} but there is
no evidence that he had discovered his famous $A$-functional
\eqref{Adef} at this stage. His audience,
especially in G\"{o}ttingen, was not always convinced that
he had sorted out all the details.
These criticisms were not put to rest until the discovery of
\eqref{Adef} and the publication of \cite{dg1931a}.
\cite{jm} contains a fuller discussion of the CCNY notebook,
the development of Douglas's ideas
and the reception he got at various seminars.

\subsection{The Riemann Mapping Theorem and
the Least Area Property} \label{rmt_lap}
An unexpected bonus of Douglas's method is a proof of
the Riemann-Carath\'{e}odory-Osgood Theorem
which follows simply by taking $n=2$. (A little work using the
argument principle is required to establish univalency of the map.)
Douglas was rightly proud that his solution not only did not require
any theorems from conformal mapping but that some such theorems could,
in fact, be proved using his method.

However, Douglas did have to use Koebe's theorem in order to establish
that his solution had least area among discs spanning $\Gamma$.
He had hoped to fix this blemish but he never succeeded.
That had to wait for contributions from Morrey \cite{morr48} and,
more recently, from Hildebrandt and von der Mosel~\cite{hild99}.
Again, details of this matter and how it contributed to the
Douglas-\rdo controversy are fully discussed in \cite{jm}.

\subsection{Solutions of the Plateau Problem of infinite area}
\label{sol_infarea} In the final section of \cite{dg1931a},
Douglas indicated that there were contours for which
every spanning surface has infinite area. Nevertheless,
he could prove the existence of a minimal surface
spanning such a contour $\Gamma$ as a limit of minimal surfaces
spanning polygonal contours which converge to $\Gamma$.
Douglas was very cross that \rdo regarded the Plateau problem
as meaningless for contours which could only bound
surfaces of infinite area. He compared the situation
to that in Dirichlet's problem, for which Hadamard had
earlier constructed continuous boundary values for which
the boundary value problem is solvable, even though
the Dirichlet functional is identically $+\infty$.
The reader is referred to \cite{jm} for a fuller account of this issue.

\subsection{Higher Connectivity and Higher Genus} Even before
working out all the details for the disc case,
Douglas was considering the Plateau problem for
surfaces of higher connectivity and higher genus.
For instance, at the February 1927 meeting of the
American Mathematical Society, \cite{dg27b}, he wrote down two integral
equations that have to be satisfied to solve the Plateau problem
for the case of two contours in $\R^n$, $n$ \emph{arbitrary.}
There is a third equation that has to be solved; it determines
the conformal type of the annulus. This is another major contribution
to the theory of minimal surfaces. As Douglas pointed out
in \cite{dg1931b}, this form of the Plateau problem
had only been raised in very special cases before
(Riemann's investigation of two parallel circles, and two
polygons in parallel planes) so his was the first general account.
It is also amusing to note that in \cite{dg1931b}, Douglas anticipated a
result reproved by Frank Morgan 50 years later in \cite{mor82}!

As early as 26 October 1929, Douglas announced that his methods
could be extended to surfaces of arbitrary genus, orientable or
not, with arbitrary many boundary curves, in a space of any
dimension.\footnote{Douglas (1930, 15, 49-50).} He may well have
had a programme at this early stage but it is doubtful that he had
complete proofs. Even when he did publish details in
\cite{dg1939}, the arguments are so cumbersome as to be
unconvincing. One should remember that Teichm\"uller theory was
still being worked out at that time and that the description of a
Riemann surface as a branched cover of the sphere is not ideally
suited for the calculation of the dependence of the
$A$-functional on the conformal moduli of the surface. Courant's
treatment in \cite{cou40} was more transparent but still awkward.
The proper context in which to study minimal surfaces of higher
connectivity and higher genus had to wait until the works of
Sacks-Uhlenbeck \cite{sau}, Schoen-Yau \cite{scy}, Jost \cite{jo}
and Tomi-Tromba \cite{tt}. Sacks-Uhlenbeck and Schoen-Yau
introduced important methods which enabled them to establish
the existence of closed minimal surfaces (without boundary)
in a closed Riemannian manifold whose universal cover
is not contractible or whose fundamental group contains
a surface subgroup. Jost extended these methods to
the boundary value problem for minimal surfaces
in Riemannian manifolds. His paper is the first published
complete solution of the Plateau-Douglas problem.
The approach of Tomi-Tromba makes use of
a differential geometric treatment of Teichm\"{u}ller space.

\section{The Fields medal}
If priority is assigned on the basis of
published papers alone, then \rdo was the first to put into print
a comprehensible solution of the Plateau problem in anything like
generality. Douglas's announcements give the impression that he
was occasionally cavalier about what he could achieve. Most of
the time he delivered on his claims, but it may not always be
appropriate to use the timing of his claims to determine priority
issues.

In \cite{jm} we conclude that \rdo and Douglas share equal credit for
solving the Plateau problem for disc-like surfaces spanning a
single contour which bounds at least one disc-like surface of
finite area. \rdo deserves full credit for solving the least area
problem. Douglas, on the other hand, was the first to solve the
Plateau problem in complete generality, that is, for an arbitrary
contour, including ones that bound only surfaces of infinite area.
He was also the only one to consider more general types of surface
than the disc, to which \rdos attentions were exclusively
confined. \rdos method for solving the Plateau problem shifts
almost all the difficulty onto problems in conformal mapping
whose solution for higher topological types was certainly not
available at the time. By contrast, Douglas's method even helped
solve some of these problems in conformal mapping. Thus,
Douglas's contributions to the Plateau problem are more major,
broader and deeper than those of Rad\'{o}. Douglas's ideas, as
developed later by Courant (who brought the Dirichlet integral
back to the forefront), have remained important in
the theory of minimal surfaces up to the present. The reader
is referred to \cite{jm} for a fuller account.

\end{document}